\documentclass[11pt, a4paper]{article}
\usepackage{latexsym}
\usepackage{indentfirst}
\usepackage{graphicx}
\usepackage{amsmath}
\usepackage{amssymb}
\begin{document}
\title{Hopf points of codimension two  in a
delay differential equation modeling leukemia}
\author{Anca-Veronica Ion$^1$, Raluca-Mihaela Georgescu$^2$\\
$^1$"Gh. Mihoc-C. Iacob" Institute of Mathematical Statistics\\
and Applied Mathematics of the Romanian Academy, Bucharest, Romania\\
$^2$University of Pite\c sti, Romania}

\date{}
\maketitle

\begin{abstract}
This paper continues the work contained in two previous papers, devoted to
the study of the dynamical system generated by a delay differential
equation that models leukemia. Here our aim  is to identify
degenerate Hopf  bifurcation points. By using an approximation of
the center manifold, we compute the first Lyapunov coefficient for
Hopf bifurcation points. We find by direct computation, in some
zones of the parameter space (of biological significance), points
where the first Lyapunov coefficient equals zero. For these  we
compute the second Lyapunov coefficient, that determines the type of the degenerate Hopf bifurcation.\\
\textbf{Acknowledgement.} Work partially supported by Grant
11/05.06.2009 within the framework of the Russian Foundation for
Basic Research - Romanian Academy collaboration.
\end{abstract}

\section{Introduction}
We consider a delay differential equation that occurs in the study of periodic
chronic myelogenous leukemia  \cite{PM-M},
\cite{PM-B-M}
\begin{equation}\label{eq}
\dot{x}(t)=-\left[\frac{\beta_0}{1+x(t)^n}+\delta\right]x(t)+k\frac{\beta_0x(t-r)}{1+x(t-r)^n}.
\end{equation}
Here  $\beta_0,\,n,\,\delta,\,k,\,r$ are positive parameters.
Parameter $k$ is of the form $k=2e^{-\gamma r}$, with $\gamma$ also
positive. We chose here to take $k$ as an independent parameter,
instead of $\gamma,$ keeping in mind the fact that, due to its
definition, $k<2$. We do not insist here on the physical
significance of the model, this is largely presented in \cite{PM-M},
\cite{PM-B-M}.

As usual in the study of differential delay equations, we consider
the Banach space $\mathcal{B}=\left\{\psi:[-r,0]\mapsto
\mathbb{R},\, \psi\, \mathrm{is\, continuous\, on\,}[-r,0]
\right\},$ with the sup norm, and for a function $x:[-r,T)\mapsto
\mathbb{R},\,T>0$ and a $0\leq t < T,$ we define the function
$x_t\in \mathcal{B}$ by $x_t(s)=x(t+s).$

In \cite{AI-stab} we proved that \eqref{eq} with the initial
condition $x_0=\phi,\,\,\phi\in \mathcal{B},$ has an unique, defined
on $[-r,\infty),$ bounded solution.

The equilibrium points of the problem are
$$x_1=0,\,\, x_2=(\frac{\beta_0}{\delta}(k-1)-1)^{1/n}.$$ The second
one is acceptable from the biological point of view if and only if
\begin{equation}\label{cond} \frac{\beta_0}{\delta}(k-1)-1
> 0,\end{equation}
condition that implies $k>1.$

In \cite{AI-stab} we studied in detail the stability of equilibrium
solutions of (\ref{eq}). The results therein are briefly recalled in Section 2.

The second equilibrium point can present
Hopf bifurcation for some points in the parameter space \cite{PM-M},
\cite{PM-B-M}, \cite{AI-Hopf}. In \cite{AI-Hopf} we considered a
typical Hopf bifurcation point,  constructed an approximation of the
center manifold in that point  and found the normal form of the
Hopf bifurcation  (by computing the first
Lyapunov coefficient, $l_1$). From this normal form we obtained the
type of stability of the periodic orbit emerged by Hopf bifurcation.
 In Section 3 we remind the ideas concerning the
approximation of the center manifold, necessary for the computation
of $l_1,$ that were developed in \cite{AI-Hopf}. In order to give a
visual representation of the Hopf bifurcation points we fix two
parameters ($n$ and $\beta_0$) and in the three dimensional
parameter space $(k,\,d,\,r)$ we represent the surface of Hopf
bifurcation points.

Let us denote by $\alpha$ the vector of parameters
$(\beta_0,\,n,\,\delta,\,k,\,r)$. If $\alpha^*$ is on the Hopf
bifurcation surface and $l_1(\alpha^*)\neq 0,$ the Hopf bifurcation
is a non-degenerate one, and $(x_2,\,\alpha^*)$ is named
\textit{codimension one Hopf point} \cite{SMB}.

If at a certain $\alpha^{**}$,  $l_1(\alpha^{**})=0,$ we have a
degenerate
 Hopf bifurcation and if we fix three of the five parameters  and vary the
other two parameters in a neighborhood of the bifurcation point in
the parameter space, we obtain a Bautin type bifurcation \cite{K}
(if the second Lyapunov coefficient, $l_2(\alpha^{**})$, is not
zero). If $l_1(\alpha^{**})=0,\,l_2(\alpha^{**})\neq 0$ then $(x_2,
\alpha^{**}),$ is also named a \emph{codimension two Hopf point}
\cite{SMB} (since Bautin bifurcation is a codimension two
bifurcation \cite{K}).

In Section 4 of the present work we present the method used by us to explore  the existence of points
with $l_1=0$ for equation \eqref{eq}. Then, for a typical such point with $l_1=0$, we develop the procedure to
compute the second Lyapunov coefficient, $l_2$. For this we compute
a higher order approximation of the center manifold.
This involves solving a set of differential equations and algebraic computation (done
 in Maple).

In Section 5 we present the results obtained for our problem by using the methods exposed in Section 4.
We found that the considered problem presents points with $l_1=0,$ and we identify
such points (to a certain approximation) by the interval bisection
technique applied with respect to one of the parameters.
We found that all Hopf codimension two points previously determined have
 $l_2<0$.

We give tables with the values of all the parameters at the points
with $l_1=0$ that we previously determined and  we  plot these
points  on the Hopf bifurcation points surface.

Section 6 presents the conclusions of our work, while Section 7 is
the Appendix containing the differential equations for the
determination of the approximation of the center manifold (and their
solutions).

\section{Stability of the equilibrium points}

The linearized equation around one of the equilibrium points is
\begin{equation}\label{lineq}
\dot{z}(t)=-[B_1+\delta]z(t)+kB_1z(t-r),\end{equation} where
$z=x-x^*,\,x^*=x_1$ or $x^*=x_2,$ $B_1=\beta'(x^*)x^*+\beta(x^*)$,
and\break $\beta(x)=\displaystyle \frac{\beta_0}{1+x^n}.$ The
nonlinear part of  equation \eqref{eq}, written in the new variable
$z$, will be denoted by $f(z(t),z(t-r)).$ The characteristic
equation corresponding to \eqref{lineq} is
\begin{equation}\label{chareq}
\lambda+\delta+B_1=kB_1e^{-\lambda r}.
\end{equation}

For the equilibrium point $x_1=0,$ we have $B_1=\beta_0$. In
\cite{AI-stab} we pointed out that for $x_1,$
$\mathrm{Re}(\lambda)<0$ for
 all eigenvalues $\lambda,$ iff the condition
$\frac{\beta_0}{\delta}(k-1)-1 < 0$ is satisfied. It follows that
$x_1$ is locally asymptotically stable in this zone of the parameter
space.

When $\frac{\beta_0}{\delta}(k-1)-1 = 0,$ $\lambda=0$ is an
eigenvalue of the linearized (around $x_1$) problem. For this case,
by constructing a Lyapunov function, we proved that $x_1$ is stable
\cite{AI-stab}. Hence, $x_1$ is stable  iff it is the only
equilibrium point (see condition \eqref{cond}). When the second
equilibrium point occurs, $x_1$ loses its stability.

For the equilibrium point $x_2,$
\begin{equation}\label{B1}
B_1=\frac{\delta}{k-1}\left[\frac{n\delta}{\beta_0(k-1)}-n+1
\right].\end{equation} The study of stability performed in
\cite{AI-stab}, that relies on the theoretical results of \cite{Ha},
reveals the following distinct situations for the stability of
$x_2$.

\textbf{I.} $B_1<0$ that is equivalent to
$\frac{\beta_0}{\delta}(k-1)>\frac{n}{n-1},$ with two subcases:

\textbf{I.A.} $B_1<0$ and $\delta+B_1<0.$ In this case $Re \lambda
<0$ for all eigenvalues $\lambda$ if and only if
$|\delta+B_1|<|kB_1|$ and
\begin{equation}\label{A}
\frac{\arccos{((\delta+B_1)/kB_1)}}{\omega_0}<r<\frac{1}{|\delta+B_1|},
\end{equation}
where $\omega_0$ is the solution in $(0,\pi/r)$ of the equation
$\omega\cot (\omega r)=-(\delta+B_1).$

\textbf{Remarks.} \textbf{1.} \textit{The condition} $\delta+B_1<0$ \textit{is equivalent to}
$\frac{\beta_0}{\delta}(k-1)(n-k)>n$ \textit{that implies} $n>k,$\textit{ and becomes}
$$\frac{\beta_0}{\delta}(k-1)>\frac{n}{n-k}. $$

\textbf{2.} \textit{The condition }$ |\delta+B_1|<|kB_1|$ \textit{is, in this case,
equivalent to} \break $\frac{\beta_0}{\delta}(k-1)>1,$ \textit{that is the condition
of existence of} $x_2,$\textit{ hence only condition }\eqref{A} \textit{brings
relevant information.}

\textbf{I.B.} $B_1<0$ and $\delta+B_1>0.$

 In this case, $Re\lambda <0$ for all eigenvalues $\lambda$ if and
only if
\begin{equation}\label{B}
\delta+B_1>|kB_1|\,\,\mathrm{or}\,\left\{\delta+B_1\leq |kB_1|\,\,
\mathrm{and}\,\,r<\frac{\arccos{((\delta+B_1)/kB_1)}}{\omega_0}\right\}
\end{equation}
where $\omega_0$ is defined as above.

Hence, when these conditions are satisfied, $x_2$ is locally
asymptotically stable.

\textbf{Remarks.} \textbf{1.} \textit{The condition} $\delta+B_1<0$ \textit{is equivalent to}
$\frac{\beta_0}{\delta}(k-1)(n-k)<n.$

\textbf{2.} \textit{Condition }$\delta+B_1>|kB_1|$ \textit{is equivalent to}
\begin{equation}\label{8}
    \frac{\beta_0}{\delta}(k-1)>\frac{n}{n-\frac{2k}{k+1}}.
\end{equation}

\textit{If }$n<\frac{2k}{k+1}$ \textit{the above inequality is satisfied since the
expression in the left hand side is positive.}

\textit{If }$n>\frac{2k}{k+1}(>1)$ \textit{then two cases are possible: either
\eqref{8} is satisfied and then }$Re \lambda < 0$ \textit{for every
eigenvalue, or}
$\frac{\beta_0}{\delta}(k-1)<\frac{n}{n-\frac{2k}{k+1}}$ \textit{and the
stability condition is satisfied if}
$r<\frac{\arccos{((\delta+B_1)/kB_1)}}{\omega_0}.$

\textbf{II.} $B_1>0$ that is equivalent to
$\frac{\beta_0}{\delta}(k-1)<\frac{n}{n-1}.$

In this case, we can only have $\delta+B_1>0,$ and, by \cite{Ha},
$Re\lambda<0$ for all eigenvalues $\lambda$ if and only if
$kB_1<\delta+B_1.$ But this inequality is equivalent to
$\frac{\beta_0(k-1)}{\delta}>1$ that is already satisfied, since
$x_2$ exists and is positive. It follows that if $B_1>0$ then $x_2$
is stable.

\section{Hopf bifurcation from the nontrivial equilibrium}

In the  stability discussion of $x_2$ shortly presented above, in
\textbf{I.}, the case
\begin{equation}\label{r-H}
r=\frac{\arccos((\delta+B_1)/(kB_1))}{\omega_0}
\end{equation}
occurs on the frontier of the stability domain.

Relations \eqref{r-H} and $\omega_0\cot(\omega_0 r)=-(\delta+B_1)$
imply
\begin{equation}\label{om0}\omega_0=\sqrt{(kB_1)^2-(\delta+B_1)^2}\end{equation} and that the pair
$\lambda_{1,2}=\pm i\omega_0$
 represents a solution of
(\ref{chareq}). All others eigenvalues have strictly negative real parts.

Thus the points in the parameter space where relations \eqref{r-H} and \eqref{om0}
hold are candidates for Hopf bifurcation. In order to have an image
of the set of such points, we fixed $n$ and $\beta_0$  and in the space $(k,\,\delta,\,r)$ we represented the surface of equation
\begin{equation}\label{surface-H}
r=\frac{\arccos((\delta+B_1)/(kB_1))}{\sqrt{(kB_1)^2-(\delta+B_1)^2}}\,.
\end{equation}
Both the numerator and the denominator of the ratio above show that the condition
$|\delta+B_1|<|kB_1|$ must be fulfilled.
In Fig. 1 we present the image of the obtained surface for $n=2,\, \beta_0=1$ and a projection of the surface
on the plane $(k,\,\delta),$ indicating the domain of definition of the function in the RHS of
\eqref{surface-H}.

\begin{figure}\centering
\includegraphics[width=0.44\linewidth]{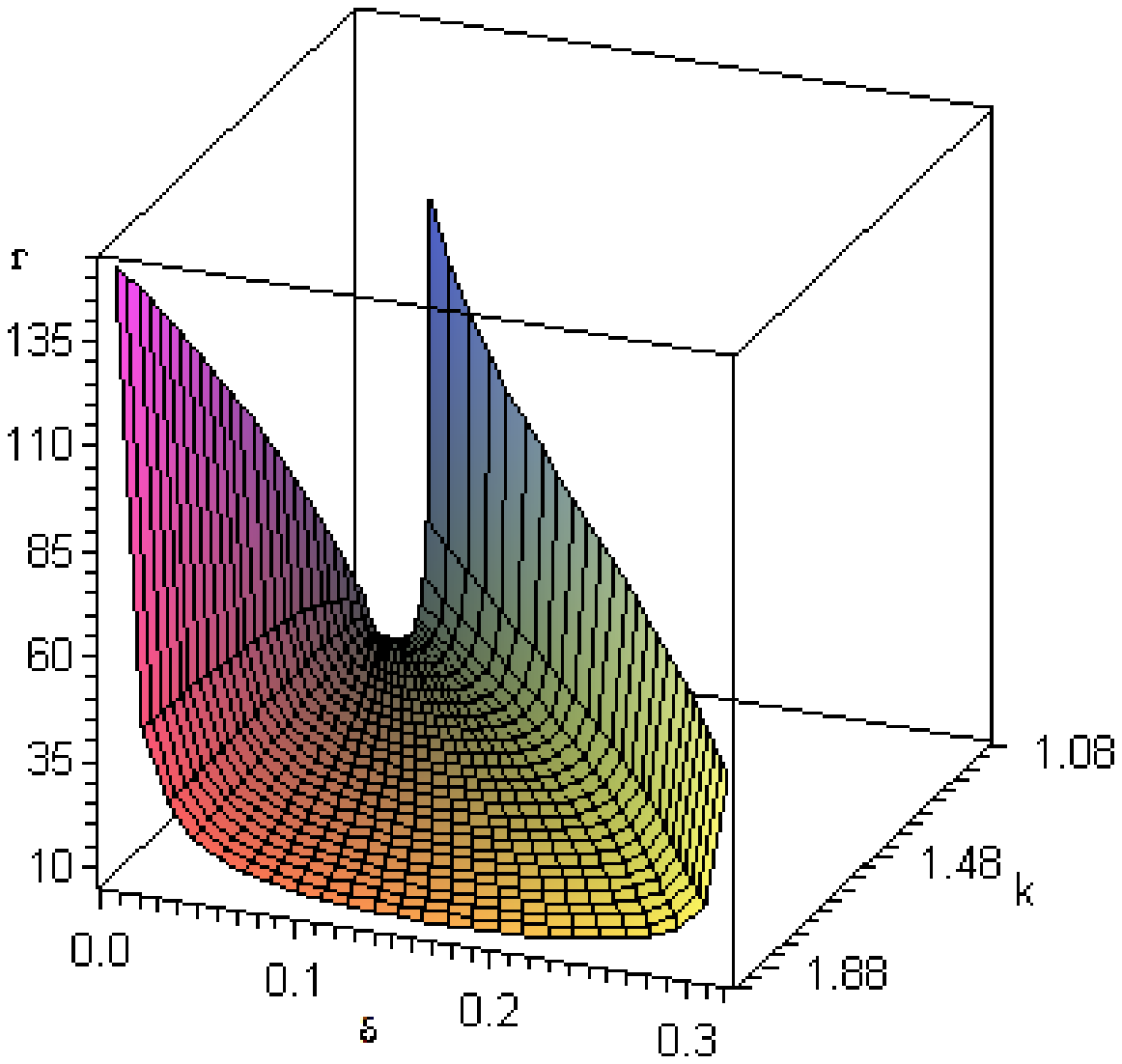}\includegraphics[width=0.49\linewidth]{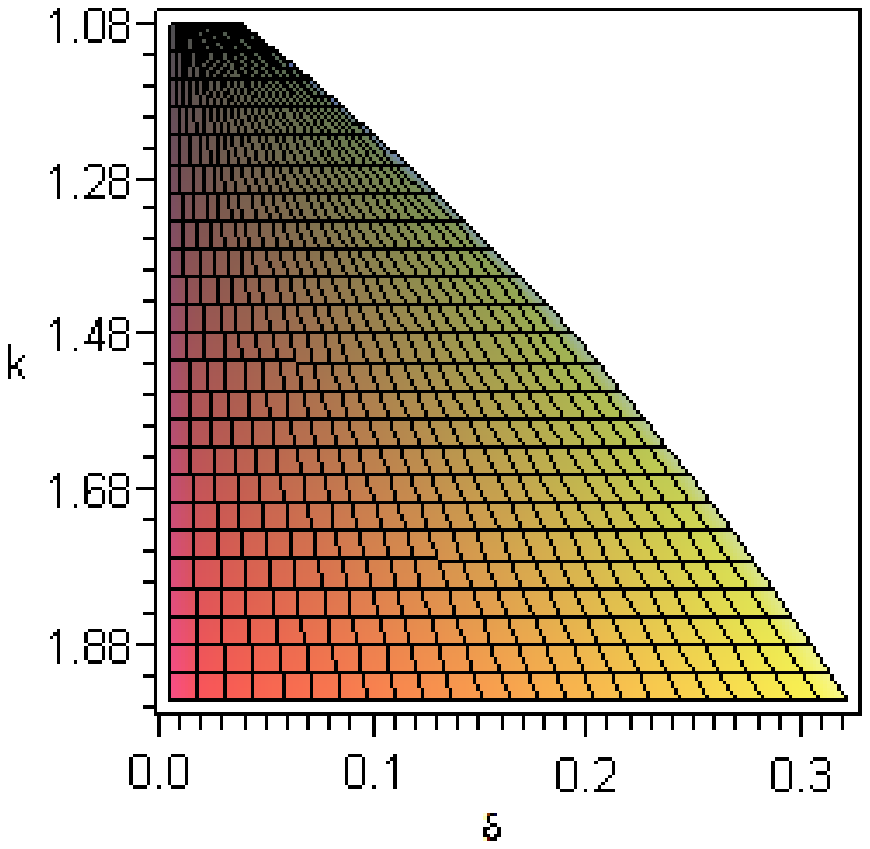}\caption{\small{\textbf{Left}- Surface of Hopf bifurcation points in the $(\delta,\,k,\,r)$ space ($n=2$,\,$\beta_0=1$). \textbf{Right}- Projection of the surface on the plane $(k,\,\delta)$ -figures done with Maple.}}
\end{figure}

Assuming that only one parameter varies (denote it by $\alpha$)
while the others are kept fixed, a point $\alpha^*,$ for which
\eqref{surface-H} is satisfied, is a nondegenerate Hopf bifurcation
point if $\frac{d\mu}{d\alpha}(\alpha^*)\neq 0$ and
$l_1(\alpha^*)\neq 0,$ where $l_1$ is the first Lyapunov
coefficient, to be defined below. If $l_1(\alpha^*)= 0,$ the Hopf
bifurcation is degenerate.

In the sequel we study the existence of points $\alpha^{**}$ in the
parameter space, satisfying \eqref{r-H} and \eqref{om0} and
$l_1(\alpha^{**})= 0,$  for the equilibrium point $x_2$. In order to
do this, we have to construct an approximation of the local
invariant center manifold and the approximation of the restriction
of problem \eqref{eq} to the local center manifold.

\subsection{The center manifold, the restriction of the problem to the center manifold}

In \cite{AI-Hopf}, by using the method of \cite{Far}, the problem
was formulated as an equation in the Banach space
$$\mathcal{B}_0=\left\{\psi:[-r,0]\mapsto \mathbb{R},\, \psi\,
\mathrm{is\, continuous\, on\,}[-r,0)\wedge\,\exists
\lim_{s\rightarrow 0}\psi(s)\in \mathbb{R}  \right\},$$ that is
\begin{equation}\label{nonlineq3}\frac{dz_t}{dt}=A(z_t)+d_0\widetilde{f}(z_t),
\end{equation}
where $\widetilde{f}(\varphi)=f(\varphi(0),\,\varphi(-r)),$
$$d_0(s)=\left\{\begin{array}{cc}
                0, & s\in[-r,0), \\
                1, & s=0, \\
              \end{array}\right.
$$
and, for $\varphi\in \mathcal{B}\subset \mathcal{B}_0,$
\begin{equation}A(\varphi)=\dot{\varphi}+d_0[L(\varphi)-\dot{\varphi}(0)],\;\;\;
L(\varphi)=-(\delta+B_1)\varphi(0)+kB_1\varphi(-r).
\end{equation}

The developments of this subsection are done for a point $\alpha^*,$ in the parameter space,  for which relations  \eqref{r-H} and \eqref{om0} are satisfied. Hence the problem has two eigenvalues $\lambda_{1,2}= \pm i \omega^*$ and all other eigenvalues have strictly negative eigenvalues.

The eigenfunctions corresponding to
 eigenvalues
$\lambda_{1,2}= \pm i \omega^*$
 are given by
$\varphi_{1,2}(s)=e^{\pm i \omega^* s},\,s\in [-r,0].$ Since the
eigenfunctions are complex functions, we need to use the
complexificates of the spaces $\mathcal{B},\,\mathcal{B}_0$ that we
denote by $\mathcal{B}_{C},\,\mathcal{B}_{0C},$ respectively. We
denote by $\mathcal{M}$ the subspace of $\mathcal{B}_C$ generated by
$\varphi_{1,2}(\cdot).$

In \cite{AI-Hopf}, by using the ideas in \cite{Far},  we constructed
a projector $\mathcal{P}:\mathcal{B}_{0C} \mapsto \mathcal{M},$ that
induces a projector defined  on $\mathcal{B}_{C}$ with values also in $\mathcal{M}$. The space $\mathcal{B}_{C}
$ is decomposed as a direct sum $\mathcal{M}\oplus \mathcal{N},$
where $\mathcal{N}=(I-\mathcal{P})\mathcal{B}_{C}.$ We do not
present here the details of this construction. The projection of problem
\eqref{eq} on $\mathcal{M}$ leads us to the equation
\begin{equation}\label{proj-eq1}
\frac{du}{dt}=i\omega^*
u+\Psi_1(0)\widetilde{f}(\varphi_1u+\varphi_2\overline{u}+\mathbf{v}),
\end{equation}
where
\begin{equation}\label{psi}\Psi_1(0)=\frac{1+(\delta+B_1-i\omega^*)r}{[1+(\delta+B_1)r]^2
+\omega^{*2}r^2}.
\end{equation}

\bigskip

The local center manifold is, for our problem, a $C^\infty$  invariant manifold,
tangent to the space $\mathcal{M}$ at the point $z=0$ (that is
$x=x_2$), and it is the graph of a function $w(\alpha^*)(\cdot)$
defined on a neighborhood of zero in $\mathcal{M}$ and taking values
in $\mathcal{N}$. A point on the local  manifold has the form
$$u\,\varphi_1+\overline{u}\,\varphi_2+w(\alpha^*)(u\,\varphi_1+\overline{u}\,\varphi_2).$$

The restriction of problem \eqref{eq} with $x_0=\phi$ to the
invariant manifold is
\begin{equation}\label{proj-eq1}
\frac{du}{dt}=\omega^*iu+\Psi_1(0)\widetilde{f}(\varphi_1u+\varphi_2\overline{u}+w(\alpha^*)(u\varphi_1+
\overline{u}\varphi_2)),
\end{equation}
with the initial condition $u(0)=u_0,$ where
$\mathcal{P}(\phi)=u_0\varphi_1+\overline{u}_0\varphi_2.$ The real
and the imaginary parts of this complex equation, represent the
two-dimensional restricted to the center manifold problem. We can
study this problem with the tools of planar dynamical systems theory
(see, e.g. \cite{K}).

\bigskip

We set
$\widetilde{w}(\alpha^*)(u,\overline{u}):=w(\alpha^*)(u(t)\varphi_1+\overline{u(t)}\varphi_2),$
and write
\begin{equation}\label{seriesw}\widetilde{w}(\alpha^*)(u,\overline{u})=
\sum_{j+k\geq2}\frac{1}{j!k!}w_{jk}(\alpha^*)u^{j}\overline{u}^{k},
\end{equation}
where
 $w_{jk}(\alpha^*)\in\mathcal{B},$  $u:[0,\infty)\mapsto
\mathbb{C}.$ By using \eqref{seriesw}, we write
\begin{equation}\label{seriesf}\widetilde{f}(\varphi_1u+\varphi_2\overline{u}+w(\alpha^*)(u\varphi_1+
\overline{u}\varphi_2))=\sum_{j+k\geq2}\frac{1}{j!k!}f_{jk}(\alpha^*)u^{j}\overline{u}^{k}.
\end{equation}

Equation  (\ref{proj-eq1}) becomes
\begin{equation}\label{eq-series} \frac{du}{dt}=\omega^*i
u+\sum_{j+k\geq2}\frac{1}{j!k!}g_{jk}(\alpha^*)u^{j}\overline{u}^{k},
\end{equation}
where
\begin{equation}\label{gij}g_{jk}(\alpha^*)=\Psi_1(0)f_{jk}(\alpha^*).\end{equation}

All the computation below are made at
$\alpha=\alpha^*$ but, for simplicity, the parameter will not be
written.

\subsection{First Lyapunov coefficient}

The first Lyapunov coefficient at $\alpha=\alpha^*$ is given by \cite{K}

\begin{equation}\label{l1}
l_1=\frac{1}{2\omega^{*2}}Re(ig_{20}g_{11}+ \omega^* g_{21}).
\end{equation}

In order to compute it,  we have to compute the corresponding
coefficients $g_{jk},$ of \eqref{eq-series}, and, for this, we must
find the coefficients of the series of $\widetilde{f}$ in
\eqref{seriesf}.

We consider the function $q(x)=\beta(x)x,$   and its derivatives of
order $n\geq 2$ in $x_2,$ i.e.
$B_n:=q^{(n)}(x_2)=\beta^{(n)}(x_2)x_2+n\beta^{(n-1)}(x_2).$ Let us
also consider the function $F(x)=B_2x^2+B_3x^3+...+B_nx^n+...$\,.

The nonlinear part  of  eq. \eqref{eq} is
$$f(z(t),z(t-r))=-F(z(t))+kF(z(t-r)),$$ or, in terms of $z_t,$ the
function $\widetilde{f}$ of (\ref{nonlineq3}) is
$$\widetilde{f}(z_t)=-F(z_t(0))+kF(z_t(-r)).$$

We denote by $\left\{\widetilde{T}(t) \right\}_{t\geq 0}$ the semigroup of operators generated by eq. \eqref{nonlineq3}. We have (by putting $\varphi_1=\varphi$, and by assuming that the initial function, $\phi,$ is on the invariant manifold)
\begin{equation}\label{termnl}\widetilde{f}(z_t)=-F([\widetilde{T}(t)\phi](0))+kF([\widetilde{T}(t)\phi](-r))=
\end{equation}
\[=-\frac{1}{2!}B_2[u\varphi(0)+\overline{u}\,\overline{\varphi}(0)+\frac{1}{2}w_{20}(0)u^{2}
+
w_{11}(0)u\overline{u}+\frac{1}{2}w_{02}(0)\overline{u}^{2}+...]^2-
\]
\[-\frac{1}{3!}B_3[u\varphi(0)+\overline{u}\,\overline{\varphi}(0)
+\frac{1}{2}w_{20}(0)u^{2}+
w_{11}(0)u\overline{u}+\frac{1}{2}w_{02}(0)\overline{u}^{2}+...]^3-....+
\]
\[+\frac{1}{2!}kB_2[u\varphi(-r)+\overline{u}\,\overline{\varphi}(-r)+
\frac{1}{2}w_{20}(-r)u^{2}+
w_{11}(-r)u\overline{u}+\frac{1}{2}w_{02}(-r)\overline{u}^{2}+...]^2+
\]
\[+\frac{1}{3!}kB_3[u\varphi(-r)+\overline{u}\,\overline{\varphi}(-r)+
\frac{1}{2}w_{20}(-r)u^{2}+
w_{11}(-r)u\overline{u}+\frac{1}{2}w_{02}(-r)\overline{u}^{2}+...]^3+...=
\]
\[=\sum_{j+k\geq 2}\frac{1}{j!k!}f_{jk}u^j\overline{u}^k.
\]
The  coefficients $f_{jk}$ involved in the expression of $l_1,$ (computed also in \cite{AI-Hopf}), are obtained by equating the same order terms in the two series from \eqref{termnl}
\[f_{20}=-B_2(1-ke^{-2\omega^* ir}),
\]
\[f_{11}=B_2(k-1),
\]
\[f_{02}=-B_2(1-ke^{2\omega^* ir}),
\]
\[f_{21}=-B_2w_{20}(0)-2B_2w_{11}(0)+2kB_2e^{-i\omega^*
r}w_{11}(-r)+kB_2e^{i\omega^* r}w_{20}(-r)-\]\[-B_3(1-ke^{-i\omega^* r}),
\]
from where $g_{ij}$ will be easily obtained from \eqref{gij} (with $\Psi_1(0)$ given in
(\ref{psi})).
Since $f_{21}$ depends on some values of $w_{20}$ and $w_{11}$, we have to determine these functions first.

The functions $w_{jk}\in \mathcal{B}_C$ are determined from the relation \cite{WWPOG},
\cite{N}, \cite{I04}
\begin{equation}\label{diffeq}\frac{\partial}{\partial s}\sum_{j+k\geq2}\frac{1}{j!k!}w_{jk}(s)u^{j}\overline{u}^{k}=
\sum_{j+k\geq2}\frac{1}{j!k!}g_{jk}u^{j}\overline{u}^{k}\varphi_1(s)+
\end{equation}
\[+
\sum_{j+k\geq2}\frac{1}{j!k!}\overline{g}_{jk}\overline{u}^{j}u^{k}\varphi_2(s)+\frac{\partial}{\partial
t}\sum_{j+k\geq2}\frac{1}{j!k!}w_{jk}(s)u^{j}\overline{u}^{k},
\]
that, by matching the same type terms yields differential equations for each $w_{jk}(s),$
while the integration constants are determined from
\begin{equation}\label{conddiffeq}
\frac{d}{dt}\sum_{j+k\geq2}\frac{1}{j!k!}w_{jk}(0)u^{j}\overline{u}^{k}+
\sum_{j+k\geq2}\frac{1}{j!k!}g_{jk}u^{j}\overline{u}^{k}\varphi_1(0)+
\sum_{j+k\geq2}\frac{1}{j!k!}\overline{g}_{jk}\overline{u}^{j}u^{k}\varphi_2(0)=
\end{equation}
\[-(B_1+\delta)\sum_{j+k\geq2}\frac{1}{j!k!}w_{jk}(0)u^{j}\overline{u}^{k}+kB_1\sum_{j+k\geq2}\frac{1}{j!k!}w_{jk}(-r)u^{j}\overline{u}^{k}+
\sum_{j+k\geq2}\frac{1}{j!k!}f_{jk}u^{j}\overline{u}^{k}.
\]

E. g. in order to determine $w_{20}(s)$ we equate the terms containing $u^2$ and obtain
\[w'_{20}(s)=2\omega^* i
w_{20}(s)+g_{20}\varphi(s)+\overline{g}_{02}\overline{\varphi}(s).\]
We integrate between $-r$ and $0$ and obtain
\[w_{20}(0)-w_{20}(-r)e^{2\omega^* ir}=
\frac{g_{20}i}{\omega^*}(1-e^{\omega^*
ir})+\overline{g}_{02}\frac{i}{3\omega^*}(1-e^{3\omega^* ir}).
\]
By equating the coefficients of  $u^2$ in \eqref{conddiffeq}, the
relation
\[\left(2\omega^* i+B_1+\delta\right) w_{20}(0)-kB_1w_{20}(-r)=f_{20}-g_{20}-\overline{g}_{02}
\]
results.
From the two relations above, we find $w_{20}(0)$ and $w_{20}(-r).$

The equations for all other functions $w_{ij}$ needed in the computation of $l_1$ and $l_2$ are all listed in the Appendix. They were obtained by symbolic computations in Maple.

\textbf{Remark.} \textit{The functions $w_{jk},\,j+k=2,$ were determined (for other
problems) by many authors \cite{SA-C}, \cite{WWPOG}, \cite{N},  in connection to the
study of Hopf bifurcations. }

Once $w_{20}(-r),\,w_{20}(0),\,w_{11}(-r),\,w_{11}(0)$ computed, we determine the value of $g_{21}$ and that of  $l_1$.

The sign of the first Lyapunov coefficient $l_1$ determines the type of the Hopf bifurcation \cite{K}. That is, if $l_1<0,$ we have a supercritical Hopf bifurcation, while if $l_1>0$ - a subcritical Hopf bifurcation.

\section{Hopf points of codimension two}

The Hopf points of codimension two are among the points of the surface of equation \eqref{surface-H}, surface that, for $n$ and $\beta_0$ fixed can be represented in $\mathbb{R}^3$ (as in Fig.1). Since they obey a new constraint, $l_1=0,$ they must lie on a curve on that surface. If we could write this constraint in a simple algebraical form, then we could attempt to represent the curve  defined by \eqref{surface-H} and $l_1=0.$ Unfortunately, this is not possible, since  the expression of $l_1,$ written in terms of the parameters of the problem is quite complex.

In this situation, our strategy for finding Hopf points of codimension two (Bautin type bifurcation points) is the
following.
\begin{itemize}
    \item As in \cite{AI-Hopf}, we chose $n^*,\,\beta^*_0,\, k^*,\,\delta^*,$
in a zone of parameters acceptable from biological point of view;
    \item with the chosen values of the parameters we compute $B^*_1$, (the value
of $B_1$ at these parameters),  next \\
                 - if $B^*_1>0,$ then $x_2$ is stable for the chosen parameters and we stop;\\
                 - if $B^*_1<0,$ we compute  $p^*=\delta^*+B^*_1,\,q^*=k^*B^*_1$;
    \item  - if $|q^*|\leq |p^*|,$ we stop;\\
           - if $|q^*|>|p^*|,$ we determine $\omega^*$ and
$r^*$ such that condition (\ref{r-H}) is fulfilled i.e. we set
$\omega^*=\sqrt{(q^*)^2-(p^*)^2},$ $r^*=\arccos(p^*/q^*)/\omega^*;$
    \item for the  found values of parameters, the first Lyapunov
coefficient, $l_1$, is computed;
    \item then we vary a parameter different of $r$,  such that the
absolute value of $l_1$ decreases - we have chosen to vary $\delta$;
    \item the above computations are repeated for several values of $\delta,$ until we
find two values $\delta_1,\,\delta_2,$ of this parameter such that
the values of $l_1$ have opposite signs - if such values exist
(obviously, at each computation, the value of $r^*$ changes, but the
condition $Re \lambda^{*}=0$ is maintained);
    \item then by using the bisection of the interval technique
(with respect to $\delta$, starting from the interval $[\delta_1,\,\delta_2]$),
we find a set of parameters such that $l_1=0$ (to a certain
approximation).
   \end{itemize}

The obtained numerical results are presented in Section 5.

\subsection{Second Lyapunov coefficient}

In order to see if the Hopf points found by the above algorithm are points of codimension two, or have a higher degeneracy, we must compute the second Lyapunov coefficient in such a point.

The formula for the second Lyapunov coefficient, at a Hopf codimension two point is \cite{K}

\[12l_2=\frac{1}{\omega^* }Re g_{32}+\frac{1}{\omega^{*2}}Im\left[g_{20}\overline{g}_{31}-
g_{11}(4g_{31}+3\overline{g}_{22})-\frac{1}{3}g_{02}(g_{40}+\overline{g}_{13})-g_{30}g_{12}\right]+
\]
\[+\frac{1}{\omega^{*3}}\left\{Re\left[g_{20}\left(\overline{g}_{11}(3g_{12}-\overline{g}_{30})+
g_{02}(\overline{g}_{12}-\frac{1}{3}g_{30})+\frac{1}{3}\overline{g}_{02}g_{03}\right)+\right.\right.
\]
\[\left.\left. + g_{11}\left(\overline{g}_{02}(\frac{5}{3}\overline{g}_{30}+3g_{12})+\frac{1}{3}g_{02}\overline{g}_{03}-
4g_{11}g_{30}\right)\right]+3Im(g_{20}g_{11})Img_{21}\right\}+\]
\[+\frac{1}{\omega^{*4}}\left\{Im\left[g_{11}\overline{g}_{02}(\overline{g}^2_{20}-3\overline{g}_{20}g_{11}-
4g_{11}^2\right]+Im(g_{20}g_{11})[3Re(g_{20}g_{11})-2|g_{02}|^2]\right\}.\]

We list below the expressions of the $f_{jk}$s corresponding to the involved $g_{jk}$s (not already given in Subsection 3.2):

\vspace{0.5cm}

\noindent$\;\,
f_{30}=-3B_2w_{20}(0)-B_{3}+3kB_2e^{-i\omega^*r}w_{20}(-r)+
kB_3e^{-3i\omega^*r}, $\\

\noindent$\;\, f_{12}=\overline{f}_{21},\,f_{03}=\overline{f}_{30},$

\[\begin{array}{ccc}
 f_{40}&=&-B_2(3w_{20}(0)^2+4w_{30}(0))-6B_3w_{20}(0)-B_4+kB_2(3w_{20}(-r)^2+  \\
 \,&\,&+4e^{-i\omega r}w_{30}(-r))+6kB_3e^{-2i\omega r}w_{20}(-r)+kB_4e^{-4i\omega r},\\
 \,&\,&\,\\
 f_{31}&=&-B_2\left(3w_{11}(0)w_{20}(0)+3w_{21}(0)+w_{30}(0)\right)+B_3(3w_{11}(0)+3w_{20}(0))-\\
 \,&\,&-B_4+kB_2(3w_{11}(-r)w_{20}(-r)+3e^{-i\omega
r}w_{21}(-r)+e^{i\omega r}w_{30}(-r))+\\
\,&\,&+kB_3(3e^{-2i\omega
r}w_{11}(-r)+3w_{20}(-r))+kB_4e^{-2\omega i r},\\
 \,&\,&\,\\
f_{22}&=&-B_2(2w_{11}(0)^2+2w_{12}(0)+w_{02}(0)w_{20}(0)+2w_{21}(0))+B_3(w_{02}(0)+                  \\
\,&\,&+4w_{11}(0)+w_{20}(0))-B_4+kB_2(2w_{11}(-r)^2+2e^{-i\omega
r}w_{12}(-r)+ \\
 \,&\,&+w_{02}(-r)w_{20}(-r)+2e^{i\omega
r}w_{21}(-r))+kB_3(w_{02}(-r)e^{-2i\omega r}+\\
 \,&\,&+4w_{11}(-r)+e^{2i\omega r}w_{20}(-r))+kB_4;
  \end{array}
\]

\noindent $\;\, f_{04}=\overline{f}_{40},\,f_{31}=\overline{f}_{31},
$
and, finally,
\[\begin{array}{ccc}
f_{32}=&-B_2\left[w_{02}(0)w_{30}(0)+6w_{11}(0)w_{21}(0)+3w_{22}(0)+3w_{12}(0)w_{20}(0)+\right.\\
 \,&\left.+2w_{31}(0)\right]-B_3\left[6w_{11}(0)^2+3w_{12}(0)+3w_{02}(0)w_{20}(0)+w_{30}(0)+\right.\\
 \,&\left.+6w_{11}(0)w_{20}(0)+6w_{21}(0)\right]+B_4\left[3w_{20}(0)+6w_{11}(0)+w_{02}(0)\right]-\\
 \,&-B_5+kB_2\left[6w_{11}(-r)w_{21}(-r)+3w_{12}(-r)w_{20}(-r)+3e^{-i\omega r}w_{22}(-r)+  \right. \\
 \,&\left.+w_{02}(-r)w_{30}(-r)+2e^{i\omega r}w_{31}(-r)\right]+kB_3\left[3w_{02}(-r)e^{-i\omega r)}w_{20}(-r)+\right.\\
 \,&+ 6e^{-i\omega r}w_{11}(-r)^2+3e^{-2i\omega r}w_{12}(-r)+e^{2i\omega r}w_{30}(-r)+\\
\,& \left.+6e^{i\omega r}w_{11}(-r)w_{20}(-r)+6w_{21}(-r)\right]+kB_4\left[3e^{i\omega r}w_{20}(-r)+\right.\\
\,&\left.+6e^{-i \omega r}w_{11}(-r)+w_{02}(-r) e^{-3i\omega
r}\right]+kB_5e^{-i\omega r}.
 \end{array}\]

The values in 0 and $-r$ of the involved $w_{ij}$s are computed by solving the differential equations and by using the additional conditions listed in the Appendix.

There is a function $w_{jk}$ that requires a special treatment, i.e.
$w_{21}$. As shown in \cite{AI-dde}, the two equations that
determine $w_{21}(0),\,w_{21}(-r)$ are dependent, and the value of
$w_{21}(0)$ may be obtained by considering a perturbed problem
(depending on a small $\epsilon$), computing the corresponding
$w_{\epsilon21}$ and by taking the limit as $\epsilon\rightarrow 0.$
The formula obtained in \cite{AI-dde} is
\begin{equation}\label{w21}w_{21}(0)=
\frac{f_{21}\langle \Psi_{1}+\Psi_{2}, \,\rho \rangle-2g_{11}\langle
\widetilde{\rho}, w_{20}\rangle -(g_{20}+2\overline{g}_{11})\langle
\widetilde{\rho}, w_{11}\rangle -\overline{g}_{02}\langle
\widetilde{\rho}, w_{02}\rangle }{2r\omega i -2rA+2},
\end{equation}
where $\langle,\,\rangle$ is the bilinear form \cite{Far},
\cite{HaL}, \cite{AI-Hopf} defined on $C([0,r],\mathbb{R})\times
C([-r,0],\mathbb{R}),$ by

\begin{equation}\label{bf}\langle\psi,\varphi\rangle=\psi(0)\varphi(0)+
kB_1\int_{-r}^0\psi(\zeta+r)\varphi(\zeta)d\zeta,
\end{equation}
$\rho(s)=-2se^{\omega is}, \;s\in [-r,0]$ and
$\widetilde{\rho}(\zeta)=-2\zeta e^{-\omega i\zeta}, \;\zeta\in
[0,r].$

After determining $w_{21}(0),$ we compute $w_{21}(-r)$ by using any
of the two relations that connect the two values. All other
necessary values of $w_{jk}$, $j+k\leq 4$ are determined by using
the relations in the Appendix.

\section{Results}

As we asserted above, we explored the equilibrium point $x_2$ with
regard to the occurrence of Hopf degenerate bifurcations. As first
step we search for points in the parameter space, where $l_1=0.$

We explored the zone of parameters $(n,\,\beta_0,\,k)$ given by:
$$(n,\,\beta_0,\,k)\in\{1,1.5,2,3,...,12\}\times\{0.5,1,1.5,2,2.5\}\times \{1.1,1.2,...,1.9\}
$$ (we remind that $k>1$ in order that $x_2$ exist, and $k<2$).

For every $(n^*,\beta^*_0,k^*)$ fixed, we looked for a $\delta^*$
and a $r^*$  such that $Re \lambda_{1,2}(\alpha^*)=0$ and $l_1=0.$

We easily see that the case $\mathbf{n=1}$ is not interesting since,
for $n=1,$  by \eqref{B1}, $B_1=\frac{\delta^2}{\beta_0(k-1)^2}>0,$
and the points with $B_1>0$ are points where $x_2$ is stable and no
Hopf bifurcations can occur.

Interesting results were found for $\mathbf{n=1.5}$ and
$\mathbf{n=2}.$ For each \break $\beta_0\in\{0.5,\,1,\,1.5,\,2,\,2.5
\},$ and each $k\in\{1.1,\,1.2,...,1.9\}$ we encountered the
following behavior: at  small values of  $\delta$ (i.e.
$\delta=0.0001$), we found $B_1<0,$ hence we could detect  Hopf
points. For these points $l_1$ was proved to be positive,  it
decreased with increasing $\delta$ and became negative for some
value of $\delta$. We then applied the bisection of interval
technique as described in Section 4.1 and actually found Hopf points
with $l_1=0$.

For each such points we computed the second Lyapunov coefficient and
we found that all points with $l_1=0,$ determined before, have
$l_2<0,$ hence they are codimension two Hopf points.

For $n=1.5$ the Hopf points occur for very large values of
$r$ ($r>60$) and, since such values for $r$ are not realistic (as
asserted in \cite{PM-M}, \cite{PM-B-M}), we present only the results
obtained for $n=2.$

In Figures 2-6  we present the sets of parameters $\alpha^{**}$ with $l_1(\alpha^{**})=0$ found by us, in tables and plotted on the surfaces of Hopf bifurcation points.
These points were determined as the points where $|l_1|$ is
approximately equal to $10^{-10}.$ A better precision may be
attained by continuing the bisection interval procedure further.

\begin{figure}[h]\centering
\begin{minipage}[h]{0.49\linewidth}
{\tiny\bf\begin{tabular}{|c|c|c|c|} \hline\
 k  & $\delta$ & $r$ & $l_2$    \\ \hline
 1.1  & 0.0045705962 & 26.125314 & -0.021 \\ \hline
 1.2  & 0.0090491351 & 25.751524 & -0.0151    \\ \hline
 1.3  & 0.0134437887 & 25.422162 & -0.0124    \\ \hline
 1.4  & 0.0177612407 & 25.130258 & -0.0108   \\ \hline
 1.5  & 0.0220070315 & 24.870352 & -0.0097   \\ \hline
 1.6  & 0.0261858065 & 24.638093 & -0.0088    \\ \hline
 1.7  & 0.0303014988 & 24.429962 & -0.0081    \\ \hline
 1.8  & 0.0343574676 & 24.243076 & -0.0076    \\ \hline
 1.9  & 0.0383566021 & 24.075039 & -0.0071    \\ \hline
\end{tabular}}
\end{minipage}
\begin{minipage}[h]{0.49\linewidth}
\includegraphics[width=\linewidth]{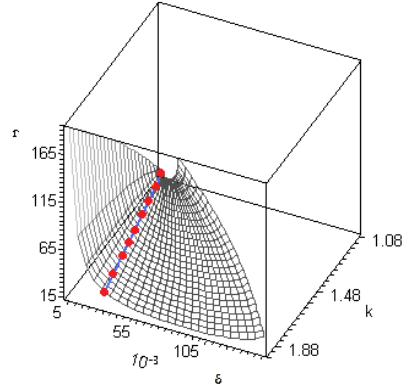}
\end{minipage}
\caption{Hopf codimension two points for
$n=2,\,\beta_0=0.5.$}\end{figure}

\begin{figure}[h]\centering
\begin{minipage}[h]{0.49\linewidth}
{\tiny\bf\begin{tabular}{|c|c|c|c|} \hline\
       k  & $\delta$ & $r$ & $l_2$    \\ \hline
 1.1 & 0.0091411924 &13.062657 & -0.0205 \\ \hline
 1.2 & 0.0180982702  &12.875762 & -0.0142\\ \hline
 1.3 & 0.0268875774 &12.711081 & -0.0114   \\ \hline
 1.4 & 0.0355224814  &12.565129 & -0.0097   \\ \hline
 1.5 & 0.0440140630  &12.435176 & -0.0085   \\ \hline
 1.6 & 0.0523716129  &12.319046 & -0.0076  \\ \hline
 1.7 & 0.0606029975  &12.214981 & -0.0069   \\ \hline
 1.8 & 0.0687149345  &12.121538 & -0.0063   \\ \hline
 1.9 & 0.0767132043  &12.037519 &  -0.0059   \\ \hline
\end{tabular}}
\end{minipage}
\begin{minipage}[h]{0.49\linewidth}\includegraphics[width=\linewidth]{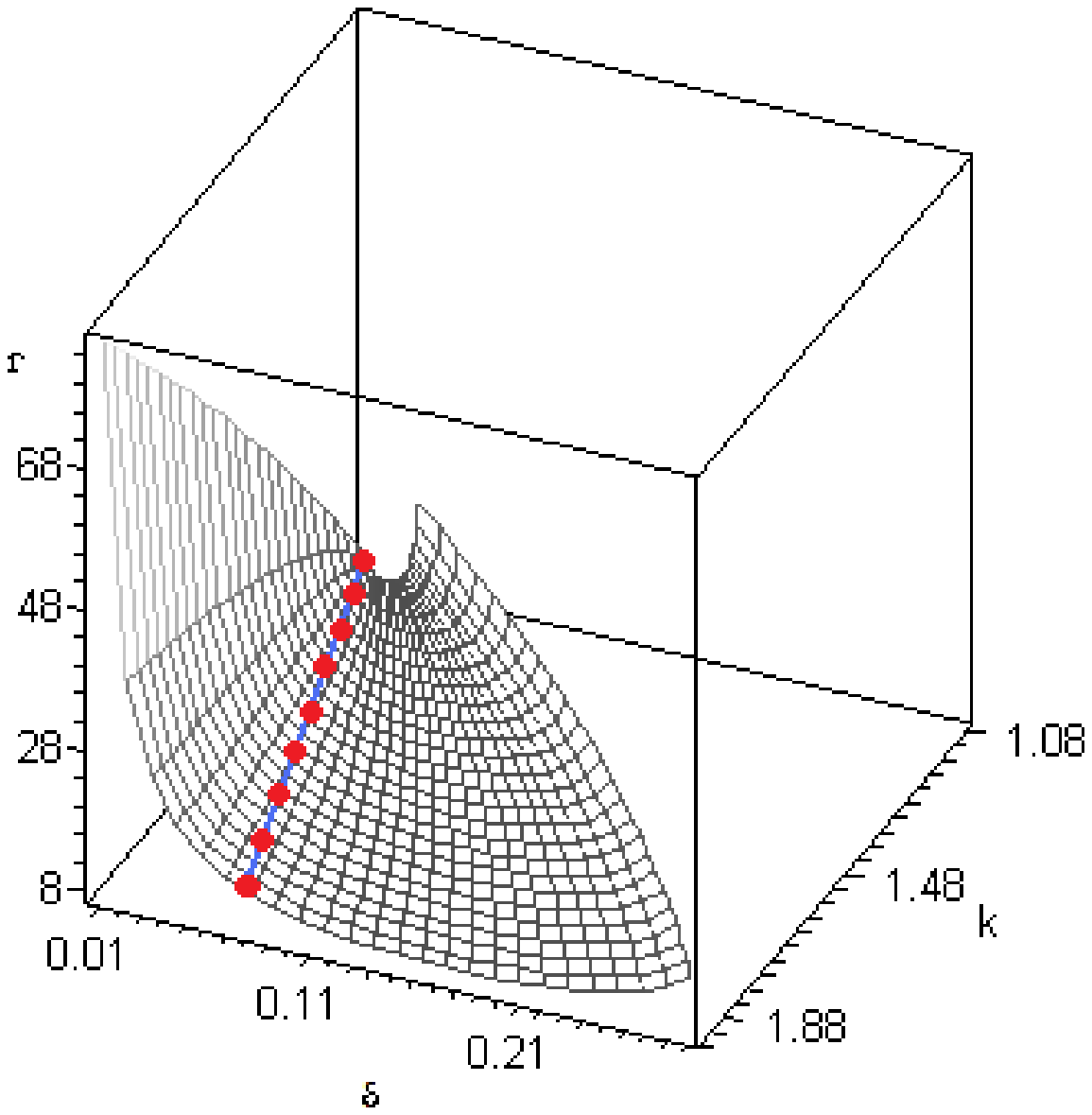}
\end{minipage}
\caption{Hopf codimension two points for
$n=2,\,\beta_0=1.$}\end{figure}

\begin{figure}[h]\centering
\begin{minipage}[h]{0.49\linewidth}
{\tiny\bf\begin{tabular}{|c|c|c|c|} \hline\
 k  & $\delta$ & $r$ & $l_2$    \\ \hline
 1.1  & 0.0137117887 & 8.708438 &  -0.0204   \\ \hline
 1.2  & 0.0271474053 & 8.583841 & -0.0140   \\ \hline
 1.3  & 0.0403313662 & 8.474054 & -0.0112   \\ \hline
 1.4  & 0.0532837222 & 8.376752 & -0.0095   \\ \hline
 1.5  & 0.0660210946 & 8.290117 &  -0.0083   \\ \hline
 1.6  & 0.0785741932 & 8.212697 & -0.0074    \\ \hline
 1.7  & 0.0909044966 & 8.143320 &  -0.0067   \\ \hline
 1.8  & 0.1030724022 & 8.081025 & -0.0061    \\ \hline
 1.9  & 0.1150698062 & 8.025013 & -0.0056 \\ \hline
\end{tabular}}\end{minipage}
\begin{minipage}[h]{0.49\linewidth}\includegraphics[width=\linewidth]{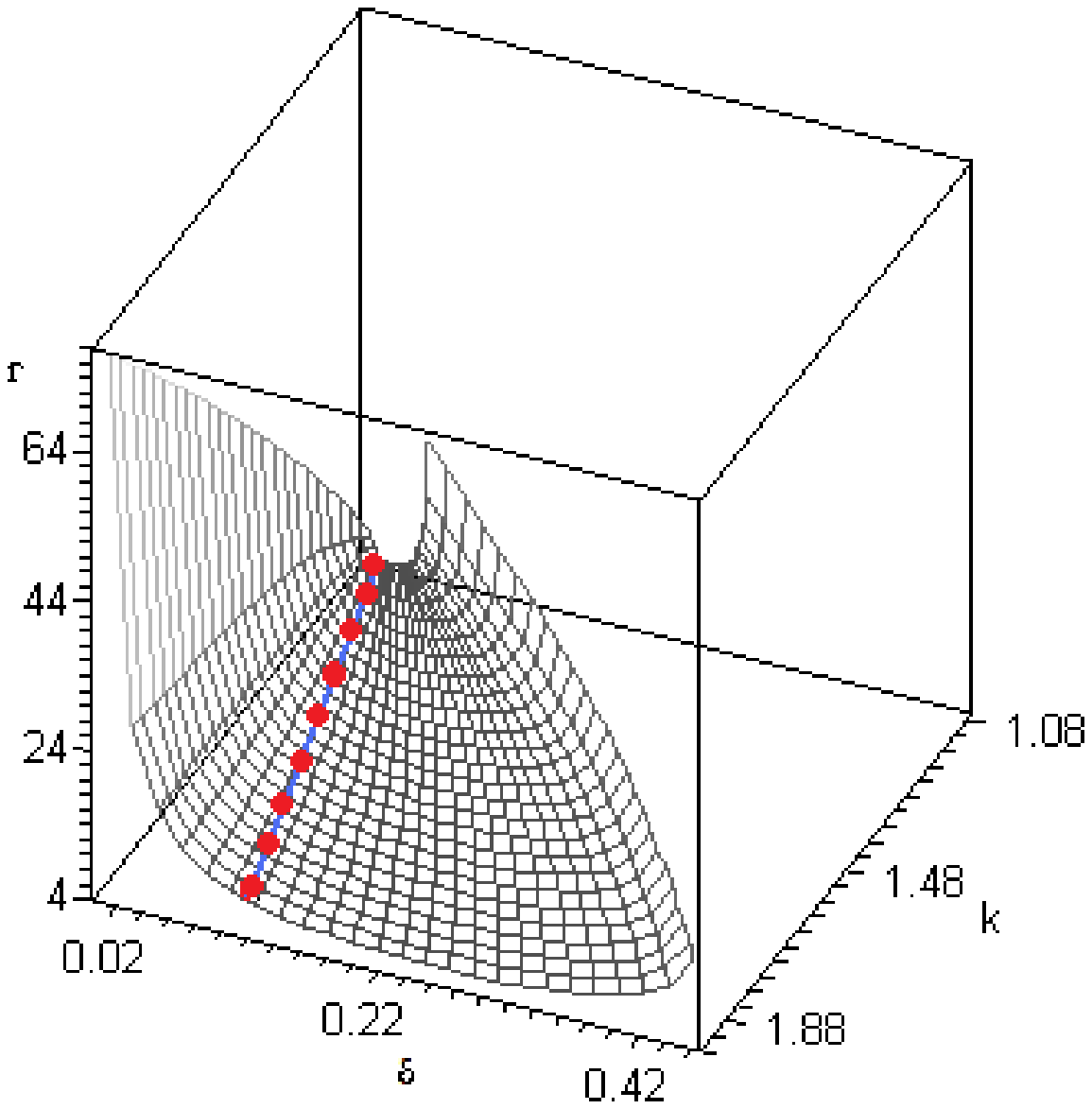}
\end{minipage}
\caption{Hopf codimension two points for
$n=2,\,\beta_0=1.5.$}\end{figure}

\begin{figure}[h]\centering
\begin{minipage}[h]{0.49\linewidth}
{\tiny\bf\begin{tabular}{|c|c|c|c|} \hline\
 k  & $\delta$ & $r$ & $l_2$    \\ \hline
 1.1  & 0.018282385 & 6.531328 & -0.0203   \\ \hline
 1.2  & 0.036196540 & 6.437880 & -0.014  \\ \hline
 1.3  & 0.053775154 & 6.355540 & -0.0111   \\ \hline
 1.4  & 0.071044963 & 6.282564 & -0.0093  \\ \hline
 1.5  & 0.088028126 & 6.217588 & -0.0082   \\ \hline
 1.6 & 0.104743225 & 6.159523 & -0.0073   \\ \hline
 1.7  & 0.121205995 & 6.107490 & -0.0066    \\ \hline
 1.8  & 0.137429869 & 6.060769 &  -0.0060   \\ \hline
 1.9  & 0.153426408 & 6.018759 & -0.0055    \\ \hline
\end{tabular}}\end{minipage}
\begin{minipage}[h]{0.49\linewidth}\includegraphics[width=\linewidth]{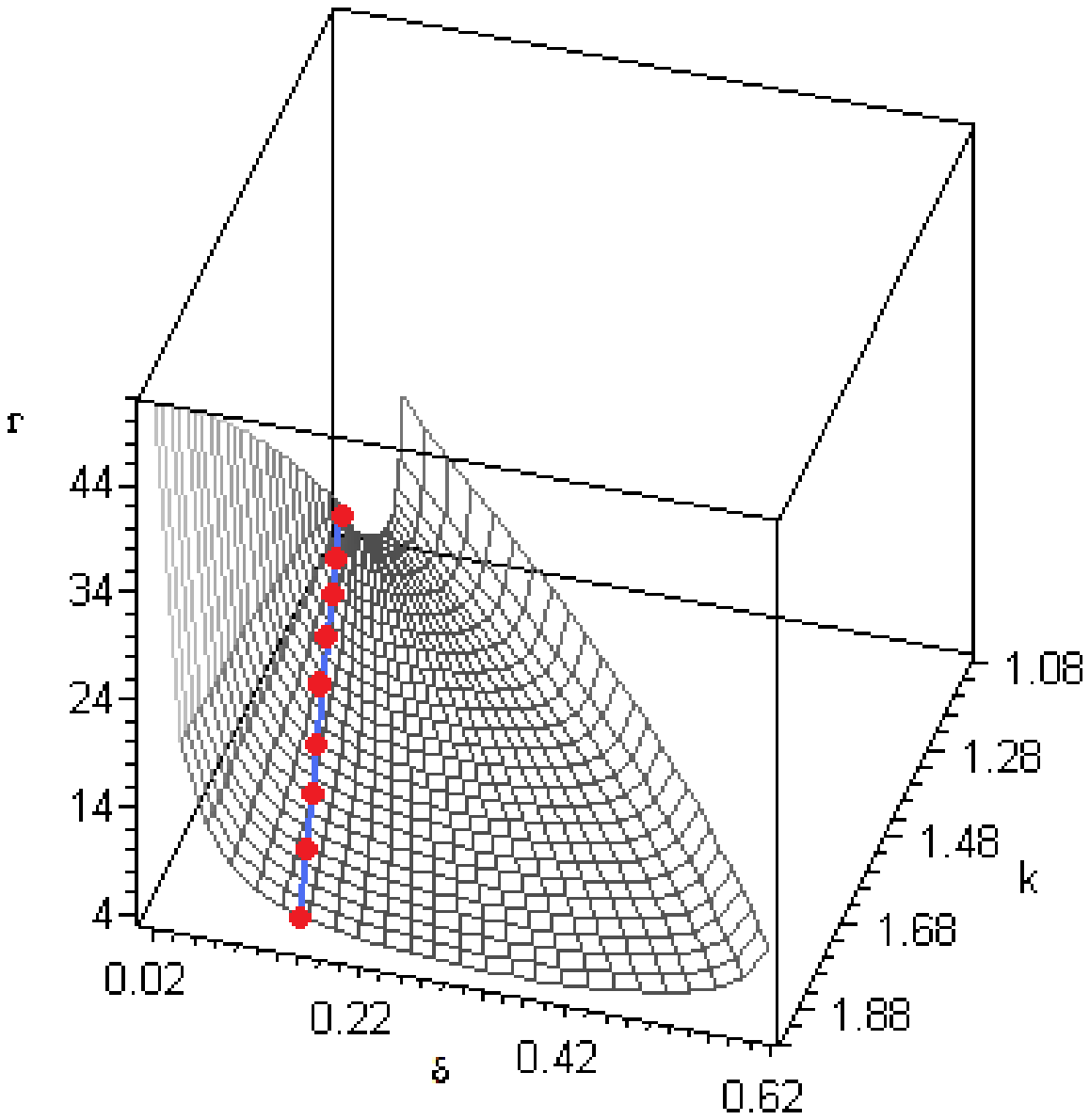}
\end{minipage}
\caption{Hopf codimension two points for
$n=2,\,\beta_0=2.$}\end{figure}

\begin{figure}[h]\centering
\begin{minipage}[h]{0.49\linewidth}
{\tiny\bf
\begin{tabular}{|c|c|c|c|} \hline\
 k  & $\delta$ & $r$ & $l_2$    \\ \hline
 1.1  & 0.022852981 &5.225062 & -0.0203  \\ \hline
 1.2  & 0.045245675 & 5.150304 & -0.0139    \\ \hline
 1.3  & 0.067218943 & 5.084432 & -0.0110    \\ \hline
 1.4  & 0.088806203 & 5.026051 & -0.0093    \\ \hline
 1.5  & 0.110035157 & 4.974074 & -0.0081    \\ \hline
 1.6  & 0.130929032 & 4.927618 & -0.0073    \\ \hline
 1.7  & 0.151507494 & 4.885992 & -0.0066    \\ \hline
 1.8  & 0.171787337 & 4.848615 & -0.0060    \\ \hline
 1.9  & 0.191783010 & 4.815007 & -0.0055    \\ \hline
\end{tabular}}
\end{minipage}
\begin{minipage}[h]{0.49\linewidth}\includegraphics[width=\linewidth]{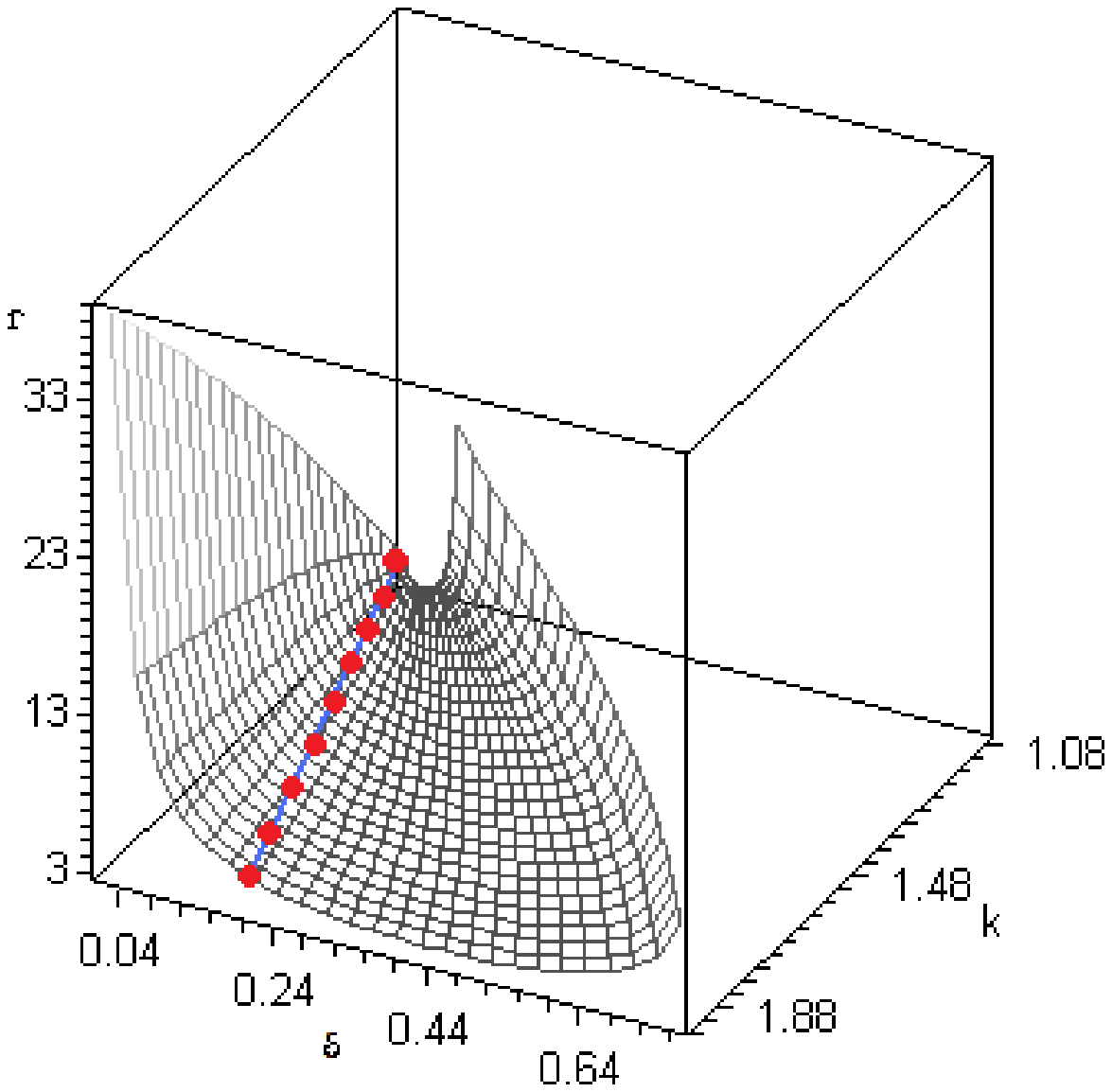}
\end{minipage}
\caption{Hopf codimension two points for
$n=2,\,\beta_0=2.5.$}\end{figure}

By connecting the points with $l_1=0$ we determine approximately a
curve of points with $l_1=0$ on each of the surfaces in Figs. 2-6.
Let us denote by $\delta=\widetilde{\delta}(k)$ the projection of
this curve on the plane $(k,\,\delta)$. Our computations of $l_1$
showed that on each of the surfaces in Figs. 2-6, $l_1>0$ for the
points of the surface that are  projected on the zone
$\delta<\widetilde{\delta}(k)$ of the plane $(k,\,\delta)$, and
hence in this zone, the Hopf bifurcation is subcritical. For the
points  that are projected on the zone
$\delta>\widetilde{\delta}(k)$ we found $l_1<0$, hence in this zone,
the Hopf bifurcation is supercritical.

Finally, for $\mathbf{3 \leq n\leq 12}$ the Hopf points found had
only negative $l_1,$ and thus, we could not find any degenerate Hopf
points.

\section{Conclusions}
The existence of Hopf codimension two points was investigated for
the non-zero equilibrium point of the equation that   models the
periodic chronic myelogenous leukemia, presented in \cite{PM-B-M},
\cite{PM-M}. We searched such points in a zone of the 5-dimensional
parameters space, of biological significance. We found that points
$\alpha$ with $l_1(\alpha)=0$ actually exist and we present tables
with values of the parameters where they occur. For each of these
points we computed the corresponding second Lyapunov coefficient
$l_2(\alpha)$ and listed the obtained values in the above mentioned
tables. We remark that for all the considered sets of parameters,
$l_2<0.$ We also remark that the values of $l_2$ are close to zero,
and, at least for $n=2$ and the chosen values of $\beta_0$ the
behavior of $l_2$ as function of $k$ is the same, i.e. $l_2$
decreases when $k$ increases. Since $k=2e^{-\gamma r},$ it can not
be greater than 2, and for $k<2$ we found only negative values of
$l_2$ (we also considered in our computation the limit case $k=2$ -
that would imply $\gamma=0$ - and for $n=2,\,k=2$ and for every
$\beta_0$ considered we found still negative $l_2$ in the points
where $l_1=0$). Hence among the considered points no higher order
degeneracies can be found (in the zones with biological
significance).

The points where $l_1=0$ and $l_2\neq 0$ are, if we vary two
parameters, points of Bautin bifurcation (provided a certain
non-degeneracy condition is satisfied \cite{I04}). A numerical
investigation of the solutions of the equation for the parameters in
and around the  points with $l_1=0,$ is in course and will be soon
submitted to publication.

\section{APPENDIX}

We list below the differential equations, their solutions, and the
supplementary conditions for the functions $w_{jk}$ with $j+k\leq4,$
excepting $w_{20}$ that is treated in Section 3. The right hand side
of the differential equations as well as the supplementary
conditions are obtained by symbolic computation in Maple.

\small{ \noindent$\mathbf{w_{11}:}$
\[\frac{d w_{11}(s)}{ds}=g_{11}e^{i\omega
s}+ \overline{g}_{11}e^{-i\omega s},
\]
\[w_{11}(s)=w_{11}(0)-\frac{ig_{11}}{\omega}(e^{\omega is}-1)+\frac{i\overline{g}_{11}}
{\omega}(e^{-\omega is}-1),
\]
\[(B_1+\delta)w_{11}(0)-kB_1w_{11}(-r)=f_{11}-g_{11}-\overline{g}_{11}.
\]
For $\mathbf{w_{02}}$, we have $w_{02}=\overline{w}_{20}$.

\noindent $\mathbf{w_{30}:}$
\[\frac{d w_{30}(s)}{ds}=3\omega i
w_{30}(s)+g_{30}e^{i\omega s}+\overline{g}_{03}e^{-i\omega
s}+3g_{20}w_{20}(s)+3\overline{g}_{02}w_{11}(s),
\]
\[w_{30}(s)=w_{30}(0)e^{3\omega is}+\frac{i}{2\omega}g_{30}(e^{i\omega
s}-e^{3\omega is})+\frac{i}{4\omega}\overline{g}_{03}(e^{-i\omega
s}-e^{3\omega is})+3g_{20}e^{3\omega is}\int_0^s
w_{20}(\theta)e^{-3\omega i\theta}d\theta+\]
\[+3\overline{g}_{02}e^{3\omega is}\int_0^sw_{11}(\theta)e^{-3\omega
i\theta}d\theta,
\]
\[(3\omega i+B_1+\delta)w_{30}(0)-kB_1w_{30}(-r)=
f_{30}-g_{30}-\overline{g}_{03}-3g_{20}w_{20}(0)-3\overline{g}_{02}w_{11}(0).
\]

\noindent $\mathbf{w_{21}:}$
\[\frac{d w_{21}(s)}{ds}=\omega i w_{21}(s)+g_{21}e^{i\omega
s}+\overline{g}_{12}e^{-i\omega s}
+2g_{11}w_{20}(s)+(g_{20}+2\overline{g}_{11})w_{11}(s)
+\overline{g}_{02}w_{02}(s),
\]
\[w_{21}(s)=w_{21}(0)e^{\omega is}+g_{21}se^{\omega is}+\frac{i}{2\omega}\overline{g}_{12}(e^{-i\omega s}-e^{\omega is})+
2g_{11}e^{\omega is}\int_0^sw_{20}(\theta)e^{-i\omega
\theta}d\theta+
\]
\[+(g_{20}+2\overline{g}_{11})e^{\omega is}\int_0^s w_{11}(\theta)e^{-i\omega \theta}d\theta
+\overline{g}_{02}e^{\omega is}\int_0^s w_{02}(\theta)e^{-i\omega
\theta}d\theta,
\]
\[(\omega i+B_1+\delta)w_{21}(0)-kB_1w_{21}(-r)=f_{21}-g_{21}-\overline{g}_{12}
-2g_{11}w_{20}(0)-(g_{20}+2\overline{g}_{11})w_{11}(0)-\overline{g}_{02}w_{02}(0).\]
\noindent $\mathbf{w_{03}:}$    $w_{03}=\overline{w}_{30}.$\\
\noindent $\mathbf{w_{12}:}$    $w_{12}=\overline{w}_{21}.$\\

\noindent $\mathbf{w_{40}:}$
\[\frac{d w_{40}(s)}{ds}=g_{40}e^{i\omega
s}+\overline{g}_{04}e^{-i\omega s}+4\omega i
w_{40}(s)+4g_{30}w_{20}(s)+4\overline{g}_{03}w_{11}(s)+6g_{20}w_{30}(s)+6\overline{g}_{02}w_{21}(s),
\]
\[w_{40}(s)=w_{40}(0)e^{4\omega
is}+\frac{i}{3}g_{40}(e^{i\omega s}-e^{4\omega
is})+\frac{i}{5}\overline{g}_{04}(e^{-i\omega s}-e^{4\omega
is})+4g_{30}e^{4\omega is}\int_{0}^{s}w_{20}(\theta)e^{-4\omega
i\theta}d\theta+
\]
\[+4\overline{g}_{03}e^{4\omega
is}\int_{0}^{s}w_{11}(\theta)e^{-4\omega
i\theta}d\theta+6g_{20}e^{4\omega
is}\int_{0}^{s}w_{30}(\theta)e^{-4\omega
i\theta}d\theta+6\overline{g}_{02}e^{4\omega
is}\int_{0}^{s}w_{21}(\theta)e^{-4\omega i\theta}d\theta,
\]
\[(4\omega i+B_1+\delta)w_{40}(0)-kB_1w_{40}(-r)=f_{40}-g_{40}-\overline{g}_{04}
-4g_{30}w_{20}(0)-4\overline{g}_{03}w_{11}(0)-\]
\[-6g_{20}w_{30}(0)-6\overline{g}_{02}w_{21}(0)
.\]

\noindent $\mathbf{w_{31}:}$

\[\frac{d w_{31}(s)}{ds}=2\omega i
w_{31}(s)+g_{31}e^{i\omega s}+\overline{g}_{13}e^{-i\omega s}+
3g_{21}w_{20}(s)+(3\overline{g}_{12}+g_{30})w_{11}(s)+
\]\[+\overline{g}_{03}w_{02}(s)+3g_{11}w_{30}(s)+3(g_{20}+\overline{g}_{11})w_{21}(s)+
3\overline{g}_{02}w_{12}(s),
\]
\[w_{31}(s)=w_{31}(0)e^{2\omega is}+\frac{i}{\omega}g_{31}(e^{i\omega
s}-e^{2\omega is})+\frac{i}{3\omega}\overline{g}_{13}(e^{-i\omega
s}-e^{2\omega is})+3g_{21}e^{2\omega
is}\int_0^sw_{20}(\theta)e^{-2\omega i\theta}d\theta+
\]
\[+
(g_{30}+3\overline{g}_{12})e^{2\omega
is}\int_{0}^{s}w_{11}(\theta)e^{-2\omega i\theta}d\theta+
\overline{g}_{03}e^{2\omega is}\int_0^sw_{02}(\theta)e^{-2\omega
i\theta}d\theta+
\]
\[+3g_{11}e^{2\omega
is}\int_0^sw_{30}(\theta)e^{-2\omega
i\theta}d\theta+(3\overline{g}_{11}+3g_{20})e^{2\omega
is}\int_0^sw_{21}(\theta)e^{-2\omega
i\theta}d\theta+3\overline{g}_{02}e^{2\omega
is}\int_0^sw_{12}(\theta)e^{-2\omega i\theta}d\theta,
\]
\[(2\omega i+B_1+\delta)
w_{31}(0)-kB_1w_{31}(-r)=f_{31}-g_{31}-\overline{g}_{13} -
3g_{21}w_{20}(0)-(g_{30}+3\overline{g}_{12})w_{11}(0)-\]
\[-
\overline{g}_{03}w_{02}(0)-3g_{11}w_{30}(0)-3(g_{20}+\overline{g}_{11})w_{21}(0)-3\overline{g}_{02}w_{12}(0).
\]

\noindent $\mathbf{w_{22}:}$
\[\frac{d w_{22}(s)}{ds}=g_{22}e^{i\omega
s}+\overline{g}_{22}e^{-i\omega
s}+2g_{12}w_{20}(s)+2(g_{21}+\overline{g}_{21})w_{11}(s) +
2\overline{g}_{12}w_{02}(s)+
\]
\[+g_{02}w_{30}(s)+(4g_{11}+\overline{g}_{20})w_{21}(s)+
(g_{20}+4\overline{g}_{11})w_{12}(s)+\overline{g}_{02}w_{03}(s),
\]
\[w_{22}(s)=w_{22}-\frac{i}{\omega}g_{22}(e^{i\omega
s}-1)+\frac{i}{\omega}\overline{g}_{22}(e^{-i\omega s}-1)
+2g_{12}\int_0^sw_{20}(\theta)d\theta+\]
\[+2(g_{21}+\overline{g}_{21})\int_0^sw_{11}(\theta)d\theta+
2\overline{g}_{12}\int_0^sw_{02}(\theta)d\theta+g_{02}\int_0^sw_{30}(\theta)d\theta+
\]
\[+(\overline{g}_{20}+4g_{11})\int_0^sw_{21}(\theta)d\theta
+(g_{20}+4\overline{g}_{11})\int_0^sw_{12}(\theta)d\theta
+\overline{g}_{02}\int_0^sw_{03}(\theta)d\theta,
\]

\[(B_1+\delta)w_{22}(0)-
kB_1w_{22}(-r)
=f_{22}-g_{22}+\overline{g}_{22}-2g_{12}w_{20}(0)-2(g_{21}+\overline{g}_{21})w_{11}(0)
-
\]
\[- 2\overline{g}_{12}w_{02}(0)-g_{02}w_{30}(0)-(4g_{11}+\overline{g}_{20})w_{21}(0)-
(g_{20}+4\overline{g}_{11})w_{12}(0)-\overline{g}_{02}w_{03}(0).
\]

\noindent $\mathbf{w_{13}:}$ $w_{13}=\overline{w}_{31}.$\\
\noindent $\mathbf{w_{04}:}$ $w_{04}=\overline{w}_{40}.$
}

\bigskip

{\small Author's addresses:\\

Anca-Veronica Ion,

"Gh. Mihoc-C. Iacob" Institute of Mathematical Statistics

 and
Applied Mathematics of the Romanian Academy,

Calea 13 Septembrie nr. 13, 050711, Bucharest, Romania.

e-mail address: anca-veronica.ion@ima.ro

\bigskip

Raluca Mihaela Georgescu,

Faculty of Mathematics and Computer Sciences,

T\^ argu din Vale, nr. 1, 110040, Pite\c sti, Arge\c s.}

 e-mail
address: gemiral@yahoo.com

\end{document}